\documentclass{amsart}

\usepackage{amsmath}
\usepackage{amsthm}
\usepackage{amssymb}

\theoremstyle{plain}
\newtheorem{theorem}{Theorem}
\newtheorem{lemma}[theorem]{Lemma}
\newtheorem{proposition}[theorem]{Proposition}
\newtheorem{corollary}[theorem]{Corollary}

\newtheorem{defn}{Definition}

\theoremstyle{remark}
\newtheorem{remark}[theorem]{Remark}

\def\gpic#1{#1
     \medskip\par\noindent{\centerline{\box\graph}} \medskip}

\DeclareMathOperator{\aut}{Aut}
\DeclareMathOperator{\lf}{leaf}

\newcommand{\ang}[1]{\langle #1 \rangle}

\author{Ian Shipman}
\title{The Distinguishing Number of the Iterated Line Graph}
\date{\today}
\email{ian.shipman@gmail.com}

\begin{document}
\begin{abstract}
We show that for all simple graphs $G$ other than the cycles $C_3,C_4,C_5,$ and the claw $K_{1,3}$ there exists a $K > 0$ such that whenever $k > K$ the $k^{th}$ iterate of the line graph can be distinguished by at most two colors.  Additionally we determine, for trees, when the distinguishing number of the line graph of $T$ is greater than the distinguishing number of $T$.
\end{abstract}

\maketitle

\section{Introduction}
Albertson and Collins introduce the concept of the distinguishing or symmetry breaking number in \cite{Alb}.  Let $G$ be a finite graph.  A map $f:V(G) \to \{1,2,\dotsc,k\}$ \emph{distinguishes} $G$ if for every nontrivial automorphism $\phi \in \aut(G)$, there exists a vertex $v \in V(G)$ where $f(v) \neq f(\phi(v))$.  The \textit{distinguishing number} of a graph, $D(G)$, is simply the least $k$ for which there exists a map $f:V(G) \to \{1,2,\dotsc,k\}$ which distinguishes $G$.  We shall refer to a map $f:V(G) \to \{1,2,\dotsc,k\}$ as a \emph{symmetry-breaking (SB)} coloring if $f$ distinguishes $G$.  Furthermore, if $k = D(G)$ then $f$ is \emph{optimal}.

For a simple graph $G$, we define the line graph $L(G)$ to be the graph whose vertices are edges of $G$ and where two edges $e,e' \in V(L(G)) = E(G)$ are adjacent if they share an endpoint in common.  We iterate the line graph in the natural way by setting $L^k(G) = L( L^{k-1}(G) )$.  When we say that a graph's iterated line graph has some property, we mean that there exists an integer $K$ such that for each $k > K$, $L^k(G)$ has the property of interest.

For a fixed group $\Gamma$, Albertson and Collins \cite{Alb} discuss the set of possible distinguishing numbers realized by a graph with automorphism group $\Gamma$.  They also show that for every finite group $\Gamma$ there is a graph $G$ with $\aut G \cong \Gamma$ and distinguishing number 2.  We note that at the heart of their construction is the the fact that $G$ behaves very rigidly under automorphisms.  $G$ can be partitioned into paths of the same length and every automorphism of $G$ acts as a bijection from the set of these paths to itself.

We will see in the third section that the iterated line graph almost always achieves the best possible distinguishing number: 2.  This result seems typical of the iterated line graph.  For instance, Knor and Niepel \cite{KnNi} show that the connectivity of the iterated line graph is optimal.  Similarly, Hartke and Higgins \cite{HH1},\cite{HH2} show that $\delta(G)$ and $\Delta(G)$ attain the minimum and maximum possible growth, respectively, in the iterated line graph.  We shall take advantage of a different sort of growth, the growth of $n(L^k(G))$, and the uniformity of this growth at each vertex.  The iterated line graph also provides a nice illustration of one interpretation of the distinguishing number.  The iterated line graph is very stiff under automorphisms, which behave by juxtaposing the elements of a partition of the vertex set and do little locally.  This property causes the iterated line graph to yield to a very simple and general method of distinguishing a graph.

In the fourth section we characterize trees whose distinguishing number increases after the first iteration. 

Additionally, the author would like extend warm thanks to Stephen Hartke for his assistance in revision, his conversations, and for suggesting the problem. 

\section{A theorem of Sabidussi}
This theorem and the map $\gamma_{G,k}$ are ultimately at the heart of our discussion.  Let $\gamma_G:\aut G \to \aut L(G)$ be given by $(\gamma_G\phi)(\{u,v\}) = \{\phi(u),\phi(v)\}$ for every $\{u,v\} \in E(G)$. We will want to extend this to arbitrary iterations of the line graph so define 
\[
	\gamma_{G,k} = \gamma_{L^{k-1}(G)} \circ \dotsm \circ \gamma_{L(G)} \circ \gamma_G
\]
which maps from $\aut G$ to $\aut L^k(G)$.  In \cite{Sab}, Sabidussi proves the following Theorem which we will use throughout.

\begin{theorem} \label{sab}
Suppose that $G$ is a connected graph that is not $P_2$, $Q$, or $L(Q)$ (see figure below).  Then $\gamma_G$ is a group isomorphism, and so $\aut G \cong \aut L(G)$.
\end{theorem}
\gpic{
\expandafter\ifx\csname graph\endcsname\relax
   \csname newbox\expandafter\endcsname\csname graph\endcsname
\fi
\ifx\graphtemp\undefined
  \csname newdimen\endcsname\graphtemp
\fi
\expandafter\setbox\csname graph\endcsname
 =\vtop{\vskip 0pt\hbox{%
    \graphtemp=.5ex
    \advance\graphtemp by 0.882in
    \rlap{\kern 0.098in\lower\graphtemp\hbox to 0pt{\hss $\bullet$\hss}}%
    \graphtemp=.5ex
    \advance\graphtemp by 0.098in
    \rlap{\kern 0.098in\lower\graphtemp\hbox to 0pt{\hss $\bullet$\hss}}%
    \graphtemp=.5ex
    \advance\graphtemp by 0.490in
    \rlap{\kern 0.882in\lower\graphtemp\hbox to 0pt{\hss $\bullet$\hss}}%
    \graphtemp=.5ex
    \advance\graphtemp by 0.490in
    \rlap{\kern 1.667in\lower\graphtemp\hbox to 0pt{\hss $\bullet$\hss}}%
    \special{pn 8}%
    \special{pa 98 882}%
    \special{pa 98 98}%
    \special{pa 882 490}%
    \special{pa 1667 490}%
    \special{fp}%
    \special{pa 98 882}%
    \special{pa 882 490}%
    \special{fp}%
    \graphtemp=.5ex
    \advance\graphtemp by 0.490in
    \rlap{\kern 2.059in\lower\graphtemp\hbox to 0pt{\hss $\bullet$\hss}}%
    \graphtemp=.5ex
    \advance\graphtemp by 0.098in
    \rlap{\kern 2.451in\lower\graphtemp\hbox to 0pt{\hss $\bullet$\hss}}%
    \graphtemp=.5ex
    \advance\graphtemp by 0.490in
    \rlap{\kern 2.843in\lower\graphtemp\hbox to 0pt{\hss $\bullet$\hss}}%
    \graphtemp=.5ex
    \advance\graphtemp by 0.882in
    \rlap{\kern 2.451in\lower\graphtemp\hbox to 0pt{\hss $\bullet$\hss}}%
    \special{pa 2059 490}%
    \special{pa 2451 98}%
    \special{pa 2843 490}%
    \special{pa 2451 882}%
    \special{pa 2059 490}%
    \special{fp}%
    \special{pa 2451 98}%
    \special{pa 2451 882}%
    \special{fp}%
    \graphtemp=.5ex
    \advance\graphtemp by 0.882in
    \rlap{\kern 1.039in\lower\graphtemp\hbox to 0pt{\hss \textit{Q} \hss}}%
    \graphtemp=.5ex
    \advance\graphtemp by 0.882in
    \rlap{\kern 3.000in\lower\graphtemp\hbox to 0pt{\hss \textit{L(Q)}\hss}}%
    \hbox{\vrule depth1.039in width0pt height 0pt}%
    \kern 3.000in
  }%
}%
}

We now have the means to gain control over the automorphisms of the iterated line graph.  Also, note that $L^2(Q)$ is not among the graphs for which Sabidussi's theorem fails to hold.  Therefore for all graphs save the paths, the sequence of automorphism groups of the iterates of the line graph eventually stabilizes.

\section{Long Term Behavior}
Suppose that $G$ is connected, simple, and not among the graphs for which Theorem \ref{sab} fails to hold.  Notice that an edge in $L(G)$ has the form $\{e,e'\}$ where $e \neq e' \in E(G)$ and $e \cap e' \neq \emptyset$.  Let us write this out using the vertices $u,v,$ and $w$.  Such an edge might look like $\{\{u,v\},\{v,w\}\}$.  Immediately, we see a natural way to associate edges in $L(G)$, i.e. vertices in $L^2(G)$, with vertices in $G$.  We define $f_G:V(L^2(G)) \to V(G)$ as
\[
f_G(\{\{u,v\},\{v,w\}\}) = v.
\]
Now if $\phi \in \aut L^2(G)$ then it follows from Theorem \ref{sab} that there exists $\phi' \in \aut G$ such that $\gamma_{G,2}(\phi') = \phi$.  Set $\ang{v}^G_1 = f_G^{-1}(v)$.

\begin{lemma} \label{sv}
Let $\phi \in \aut L^2(G)$ and $\phi' \in \aut G$ such that $\gamma_{G,2}(\phi') = \phi$.  Then $\phi(\ang{v}^G_1) = \ang{\phi'(v)}^G_1$ for all $v \in V(G)$.
\end{lemma}
\begin{proof}
Let $\psi = \gamma_G(\phi')$.  Then
\begin{multline*}
 \phi(\{u,v\},\{v,w\}) = \{\psi(\{u,v\}),\psi(\{v,w\})\} \\
 = \{\{\phi'(u),\phi'(v)\},\{\phi'(v),\phi'(w)\}\} \in \ang{\phi'(v)}^G_1.
\end{multline*}
Similarly if $\{\{u,\phi'(v)\},\{\phi'(v),w\}\} \in \ang{\phi'(v)}^G_1$ then take $u' = \phi'^{-1}(u)$ and $w' = \phi'^{-1}(w)$ to get
\begin{multline*}
	\{\{u,\phi'(v)\},\{\phi'(v),w\}\} = \{\{\phi'(u'),\phi'(v)\},\{\phi'(v),\phi'(w)\}\} \\
	= \phi( \{\{u',v\},\{v,w'\}\}) \in \phi(\ang{v}^G_1).
\end{multline*}
\end{proof}

Let us pull these ideas out past 2 iterations.
\begin{defn}
Let $\ang{v}^G_1 = f_G^{-1}(v)$ and
\[
\ang{v}^G_{k+1} = \bigcup_{u \in \ang{v}^G_{k}}{f_{L^{2k}(G)}^{-1}(u)}.
\]
If the graph $G$ is clear from context we shall denote $\ang{v}^G_k$ by $\ang{v}_k$.
\end{defn}

\begin{corollary}
Let $m >0$, $\phi \in \aut L^{2m}(G)$ and $\phi' \in \aut G$ such that $\gamma_{G,2m}(\phi') = \phi$.  Then $\phi(\ang{v}_m) = \ang{\phi'(v)}_m$ for all $v \in V(G)$.
\end{corollary}
\begin{proof}
We use induction on $m$.  When $m = 1$, this is Lemma \ref{sv}.  Suppose that $\phi$ and $\phi'$ are as in the hypotheses.  Let $\psi = \gamma_{G,2m-2}(\phi')$.  Then by the inductive hypothesis, $\psi( \ang{v}_{m-1} ) = \ang{\phi'(v)}_{m-1}$.  By Lemma \ref{sv} we also know that $\phi(f^{-1}_{L^{2m-2}(G)}(v)) = f^{-1}_{L^{2m-2}(G)}(\psi(v))$ for each $v \in V(L^{2m-2}(G))$.  Hence we have
\[
	\phi(f^{-1}_{L^{2m-2}(G)}(\ang{v}_{m-1})) = f^{-1}_{L^{2m-2}(G)}(\psi(\ang{v}_{m-1}))= f^{-1}_{L^{2m-2}(G)}(\ang{\phi'(v)}_{m-1}).
\]
But $f^{-1}_{L^{2m-2}(G)}(\ang{v}_{m-1}) = \ang{v}_m$ and $f^{-1}_{L^{2m-2}(G)}(\ang{\phi'(v)}_{m-1}) = \ang{\phi'(v)}_m$, yielding the desired equality.
\end{proof}

\begin{lemma}
If $\delta(G) \geq 3$ then for all $v \in V(G)$ and $m \geq 1$ we have $|\ang{v}_{m+1}| > |\ang{v}_m|$.
\end{lemma}
\begin{proof}
We recall that if $e = \{a,b\} \in E(G)$ then $deg_{L(G)}(e) = deg_G(a) + deg_G(b) - 2$.  Hence $\delta(L(G)) \geq 2\delta(G) - 2 \geq 4$ and repetition of this argument yields $\delta(L^{2m}(G)) \geq 3$.  Since $|N(z)| \geq 3$ for each $z \in V(L^{2m}(G))$, there exist $u,w,x \in N(z)$ and therefore we can immediately produce
\begin{align*}
\{\{u,z\},\{z,w\}\} \\
\{\{u,z\},\{z,x\}\} \\
\{\{x,z\},\{z,w\}\}
\end{align*}
in $f^{-1}_{L^{2m}(G)}(z)$.  Now we recall that
\[
	\ang{v}_{m+1} = \bigcup_{u \in \ang{v}_m}{f_{L^{2m}(G)}^{-1}(u)},
\]
and because $\delta(L^{2m}(G)) \geq 3$ we have $f^{-1}_{L^{2m}(G)}(u) \geq 3$ for all $u \in V(L^{2m}(G))$.  Therefore $|\ang{v}_{m+1}| > |\ang{v}_m|$ as desired.
\end{proof}

If our graph is suitably free of tendrils, then iterating the line graph causes each vertex $v$ to blossom into a cluster of vertices $\ang{v}_m$.  In fact, if $\delta(G) \geq 3$, there are at least $3^m$ vertices associated with $v$ in $L^{2m}(G)$.  However, the automorphisms of the iterated line graph must map clusters associated with one vertex onto clusters associated with another.  Thus, with respect to symmetry, the graph behaves as though the clusters were identified with their ancestor vertices.  Once these clusters get large enough, we can simply number them uniquely which effectively kills the nontrivial automorphisms.  We need one more easy Lemma though.  (For one proof, see \cite{KNS}.)

\begin{lemma} \label{mindegl}
Suppose that $G$ is a connected graph which is not a path, cycle or $K_{1,3}$.  Then for some $k$, $\delta(L^k(G)) \geq 3$.
\end{lemma}

\begin{theorem} \label{mainthm}
If $G$ is a connected graph other than $C_3,C_4,C_5$ or $K_{1,3}$ then there exists a $K$ such that for all $k \geq K$ we have $D(L^k(G)) \leq 2$.
\end{theorem}
\begin{proof}
First, note that if $G$ is a cycle of length at least 6 or a path then $D(G) = 2$ already.  The line graph operator fixes cycles so that $D(L^k(C_n))=2$ for all $k \geq 0$ and all $n \geq 6$.  Also, the line graph operator takes paths to paths and therefore $D(L^k(P_n)) =2$ for all $k < n$ and $D(L^k(P_n))=0$ for all $k \geq n$.  Now, we assume that $G$ is not a path, cycle, or $K_{1,3}$.  Let $H = L^m(G)$ be such that $\delta(H) \geq 3$ and $\aut(H) = \aut(L(H))$.  Such an $m$ exists by Lemma \ref{mindegl} and Theorem \ref{sab}.  Then for all $r\geq 1$ and $v \in V(H)$, $|\ang{v}^H_{r+1}| > |\ang{v}^H_r|$.  Choose $r$ such that $\min_{v \in V(H)}{|\ang{v}^H_r|} \geq n(H)$.  Then we color $V(L^{2r}(H))$ as follows.  Let $\chi:V(L^{2r}(H)) \to \{0,1\}$ be such that 
\[
	\sum_{w \in \ang{u}^H_r}{\chi(w)} \neq \sum_{w \in \ang{v}^H_r}{\chi(w)}
\]
for all $u \neq v \in V(H)$.  This simply requires using a different number of ones to color each of the $\ang{v}_r$ and we are guaranteed that this is possible because (1) for all $v \in V(H)$, $|\ang{v}^H_r| > n(H)$ and (2) for $u \neq v \in V(H)$, $\ang{u}^H_r \cap \ang{v}^H_r = \emptyset$.  Suppose that $\phi \neq id \in \aut L^{2r}(H)$ and let $\phi' = \gamma_{H,2r}^{-1}(\phi)$.  There exists $v$ such that $\phi'(v) \neq v$.  But $\phi(\ang{v}^H_r) = \ang{\phi'(v)}^H_r$ and $\ang{v}^H_r \neq \ang{\phi'(v)}^H_r$. Since
\[
	\sum_{w \in \ang{v}^H_r}{\chi(w)} \neq \sum_{w \in \ang{\phi'(v)}^H_r}{\chi(w)}
\]
the coloring we presented cannot possibly agree with $\phi$.  Of course, $\phi$ was arbitrary; hence $\chi$ distinguishes $L^{2r}(H) = L^{m + 2r}(G)$.  Clearly we can repeat this coloring procedure for $L^{m + 2s}(G)$ where $s \geq r$.  However, we also note that since $\delta(L(H)) \geq 3$ we can also find an $r'$ such that $L^{(m+1) + 2r'}(G)$ has distinguishing number at most 2 and the same is true.  Therefore we set $K = \max\{m + 2r, (m+1) + 2r'\}$ so that for each $k \geq K$ we get $D(L^k(G)) \leq 2$.
\end{proof}

\begin{remark}
Often, we can improve $K$ by choosing $r$ such that $\min_{v \in V(H)}{|S_{v,r}|} \geq D(H)$.  We modify the coloring scheme by first choosing an optimal SB coloring $f$ of $H$.  Then we color $L^{2r}(H)$ so that $\sum_{w \in S_{u,r}}{\chi(w)} = f(u)$.  Since $D(H) \leq n(H)$ we stand to get a smaller $r$ and thus a smaller $k = m + 2r$.  The same principle applies to odd iterates of the line graph.  Now, consider the subgraph $H$ of $G$ obtained by deleting all the vertices of degree at least 3.  Let $p$ be the size of the largest component.  Then we easily bound 
\[
	K \leq \max\{p + 2\log_3(D(L^p(G))), (p + 1) + 2\log_3(D(L^{p+1}(G)))\}. 
\]

\end{remark}

We now move on to the analysis of trees under one iteration of $L$, to get a taste of what can happen.

\section{Short Term Behavior: Trees}

Trees have very easy to understand automorphism groups.  For this reason, we can learn quite a bit about how the task of distinguishing a tree changes when we move to its line graph.  This close up analysis provides a nice contrast with the result of the second section.  For convenience we will fix one tree and proceed to analyze it.  Let $T$ be this tree and let $C = C(T)$ denote the center of $T$.  The center will play a large role in our discussion because it acts as a ``pivot point'' in the sense that it is fixed (as a set) under all the automorphisms of the tree.  Trees may be seen as a collection of ``branches'' attached to this central pivot point.

Toward developing an understanding of tree automorphisms, let $\lf(T)$ be the set of leaves of $T$.  Now we can define a sequence of subtrees $T = T_0, T_1,T_2,\dotso$ by letting $T_{i+1} = T_i \setminus \lf(T_i)$.  So each new subtree comes from the previous subtree by removing its leaves.  The automorphisms of $T$ fix these subtrees as sets of vertices.

\begin{lemma} \label{leaves}
For any $\phi \in \aut T$ we have $\phi(\lf(T_i)) = \lf(T_i)$.
\end{lemma}
\begin{proof}
We use induction on $i$.  Since $deg(v) = deg(\phi(v))$ we must have $\phi(\lf(T_0)) = \lf(T_0)$.  Now consider $v \in \lf(T_i)$.  We note that since $\phi(\lf(T_j)) = \lf(T_j)$ for all $i < j$, we have $\phi(T_i) = T_i$ and therefore since $\phi|_{T_i}$ is an automorphism, the degree requirement again requires $\phi(v) \in \lf(T_i)$.
\end{proof}

Notice that if $v \in V(T)$ is a leaf then there is only one edge incident on $v$ and therefore there is a correspondence between leaves of $T$ and certain vertices in $L(T)$.  We can extend this correspondence into the tree.  For each $T_i$ we associate the leaves of $T_i$ with their edges if any edges are present.  Now if $T$ has one vertex as the center then this center vertex does not have an edge associated with it and if $T$ has an edge as the center then both endpoints of this edge are associated with the same edge.  For every vertex $v \in V(T)$ that is not the center, let $e_v$ be the edge associated with $v$.  This correspondence is always onto $E(T)$.  We note that if $e_v = \{v,u\}$ then for some $i$, $v \in \lf(T_i)$ and $u \in \lf(T_{i+1})$ so if $e=\{u,v\} \in E(T)$ and $v \in \lf(T_i)$ while $u \in \lf(T_{i+1})$ then $e$ is associated with $v$ and not with $u$.

\begin{lemma}
If $\phi' \in \aut L(T)$ and $\phi \in \aut T$ is such that $\gamma_T(\phi) = \phi'$ then $\phi'(e_v) = e_{\phi(v)}$ for all $e_v \in E(T)$.
\end{lemma}
\begin{proof}
Let $e_v = \{v,u\}$.  Then $\phi'(e_v) = \{\phi(v),\phi(u)\}$.  We consider two cases.  If $e_v = C(T)$ then $e_v = e_u = e_{\phi(v)} = e_{\phi(u)}$.  If $e_v \neq C(T)$ then $v \in \lf(T_i)$, $u \in \lf(T_{i+1})$ and since $\phi(\lf(T_i)) = \lf(T_i)$ we have $\phi(v) \in \lf(T_i)$ and $\phi(u) \in \lf(T_i)$ and therefore $\{\phi(v),\phi(u)\} = e_{\phi(v)}$.
\end{proof} 

\begin{proposition}
For any tree $T$ other than $P_2$, $D(T) \leq D(L(T))$.
\end{proposition}
\begin{proof}
First, if $D(L(T)) = 1$ then $\aut(L(T))$ is trivial and since $T \neq P_2$ this implies that $\aut(T)$ is also trivial and therefore that $D(T) = 1$ as well.  So suppose that $k \geq 2$ and $\chi:V(L(T)) \to [k]$ is an optimal SB coloring of $L(T)$.  Define $\chi':V(T) \to [k]$ by the rule $\chi'(v) = \chi(e_v)$.  We then modify $\chi'$ depending on the center of $T$.  If the center of $T$ is a vertex we just color it 0.  If the center of $T$ is an edge we color its endpoints 0 and 1.  We claim that $\chi'$ is an SB coloring for $T$.  Suppose that $\phi \in \aut T$ is not the identity.  Then $\phi' = \gamma_T(\phi)$ is not the identity in $\aut L(T)$.  So there exists an edge $e = e_v$ such that $e_v \neq \phi'(e_v) = e_{\phi(v)}$ and $\chi(e_v) \neq \chi(e_{\phi(v)})$.  But $\chi(v) = \chi'(e_v) \neq \chi'(e_{\phi(v)}) = \chi(\phi(v))$.  Since $\phi$ was arbitrary, the only automorphism in $\aut T$ that agrees with $\chi$ is the identity.
\end{proof}

\begin{proposition} \label{monocent}
If $T$ is a tree which has an optimal SB coloring in which the center receives a single color then $D(T) = D(L(T))$.
\end{proposition}
\begin{proof}
Suppose that $\chi:V(T) \to [k]$ is an optimal SB coloring in which the center receives a single color.  We define an SB coloring of $L(T)$ by $\chi':V(L(T)) \to [k]$ by the rule $\chi'(e_v) = \chi(v)$, and this is unambiguously defined since the center is monochromatic.  Hence $D(T) \geq D(L(T))$ and together with the previous lemma this implies $D(T) = D(L(T))$.
\end{proof}

Above, we ask that the tree effectively has a trivial center, by requiring it to be one-colorable under an optimal SB coloring.  For $n > 1$ we can construct a tree $T$ as follows.  Begin with an edge.  To each of the endpoints append $n$ leaves.  This tree has no optimal SB coloring in which the center is monochromatic and $D(L(T)) = D(T) + 1$.  This construction can be generalized to a large family of ``saturated'' trees using the ideas in proposition \ref{treechar}.

Let $S = \{S_1,S_2,\dotsc,S_k\}$ be the collection of components of $T - C$.  Let $v_1,v_2,\dotsc,v_k$ be the vertices $v_i \in V(S_i)$ which are adjacent to vertices in $C$.  We define an equivalence relation on $S$ by $S_i \sim S_j$ when there exists an isomorphism $\beta_{i,j}:S_i \to S_j$ such that $\beta_{i,j}(v_i) = v_j$.  Let $P = \{P_1,P_2,\dotsc,P_r\}$ be the set of equivalence classes under $\sim$.  

We now consider a concept related to distinguishing.  Let $S_i \in S$ and consider the subgroup $A_i \leq \aut S_i$ which fixes $v_i$.  Then there is a distinguishing number associated with $A_i$ as well.  However, we are more interested in the \textit{number} of colorings that distinguish $S_i$ with respect to the automorphisms in $A_i$.  Suppose that $f_1,f_2$ are two $k$-colorings of $S_i$.  Then we say that $f_1$ and $f_2$ are equivalent when there is an automorphism $\phi \in A_i$ such that $f_1(v) = f_2(\phi(v))$ for all $v \in V(S_i)$.  Consider the set of $k$-colorings of $S_i$ that distinguish $S_i$ with respect to $A_i$.  This notion of equivalence produces equivalence classes in this set.  We define $D(S_i;k)$ to be the number of such equivalence classes.  An equivalence class is then a rooted SB class.  Then $D(S_i;k)$ is the number of essentially distinct $k$-colorings which distinguish $S_i$ with respect to $A_i$. 

In preparation for our next theorem, let us develop an idea of how features of $T$ persist into $L(T)$.  Let $S_i \in S$.  Then we define $S'_i$ to be the subgraph of $L(T)$ with vertex set $\{e_v: v \in V(S_i)\}$ and all the edges from $L(T)$ which connect pairs of vertices in this set.  We simply associate $S_i$ with its edges and add the edge which connects $S_i$ to $C$ so that $n(S_i) = n(S'_i)$.  Let $A'_i \leq \aut S'_i$ be the collection of automorphisms which fix $e_{v_i}$.  Note that $A'_i \cong A_i$.  We define $D(S'_i;k)$ analogously to $D(S_i;k)$ save that $e_{v_i}$ takes the place of $v_i$ as the fixed point.  If $s = S_j \in S$ then for convenience we refer to $v_j$ as $v_s$. 

\begin{lemma} \label{presdprime}
Suppose that $v_i$ is adjacent to $u \in C(T)$.  Then $D(S_i;k) = D(S'_i;k)$.
\end{lemma}
\begin{proof}
We use a modified version of the coloring procedure of proposition \ref{monocent}.  Suppose that $f$ is a k-coloring of $S_i$ that distinguishes $S_i$ with respect to $A_i$.  Define a $k$-coloring $f'$ by $f'(e_v) = f(v)$.  Suppose that $\phi' \in A'_i$ be a non-identity automorphism and let $\phi \in A_i$ be such that $\gamma_{S_i}(\phi) = \phi'$.  Then there exists $v \in V(S_i)$ such that $f(v) \neq f(\phi(v))$.  Then $f'(e_v) = f(v) \neq f(\phi(v)) = f'(e_{\phi(v)})$.  But $\phi'$ was arbitrary so this coloring distinguishes $S'_i$.  Clearly if $f,g$ are two different colorings then $f'$ and $g'$ will be different.  Hence, $D(S_i;k) \leq D(S'_i;k)$ and other direction is nearly identical.
\end{proof}

Prior to the next theorem, observe that if $A,B \in P_i$ then $D(A;k) = D(B;k)$ since $A$ and $B$ are isomorphic as rooted trees.

\begin{proposition}  \label{treechar}
Let $T$ be a tree and $k = D(T)$.  Using the definitions above, for each $P_i$ set $m_i = D(S_j;k)$ where $S_j \in P_i$.  Suppose $C(T) = \{u,w\}$ then define $U,W \subset S$ to be those components of $T - C(T)$ which are adjacent to $u$ and $w$ respectively.  Then $D(T) < D(L(T))$ if and only if for all $i$, $|P_i \cap U| = |P_i \cap W| = m_i$.
\end{proposition}

\begin{remark}
When the hypotheses of this proposition are met, the tree is as full on either side of the center as can be and only by distinguishing the two vertices in the center can we attain the distinguishing number given.  In the line graph the center collapses to a point and we can map one side of the line tree to the other since all possible distinct rooted SB $k$-colorings on subtrees were used on both sides.  As we have seen if the center of $T$ is an edge, then the center collapses to a point.  The proof idea below extends naturally to a similar statement about $T \cdot e$, the contraction of $e$, where $e$ is the center of $T$, namely $D(T) \leq D(T\cdot e)$.
\end{remark}

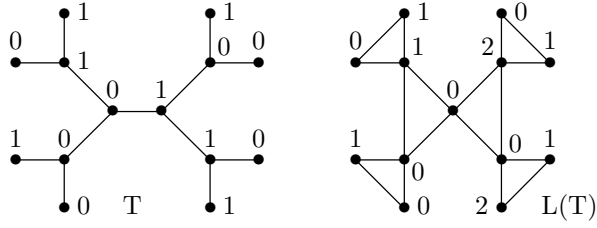
\gpic{
\expandafter\ifx\csname graph\endcsname\relax
   \csname newbox\expandafter\endcsname\csname graph\endcsname
\fi
\ifx\graphtemp\undefined
  \csname newdimen\endcsname\graphtemp
\fi
\expandafter\setbox\csname graph\endcsname
 =\vtop{\vskip 0pt\hbox{%
    \graphtemp=.5ex
    \advance\graphtemp by 1.119in
    \rlap{\kern 0.356in\lower\graphtemp\hbox to 0pt{\hss $\bullet$\hss}}%
    \graphtemp=.5ex
    \advance\graphtemp by 0.864in
    \rlap{\kern 0.356in\lower\graphtemp\hbox to 0pt{\hss $\bullet$\hss}}%
    \graphtemp=.5ex
    \advance\graphtemp by 0.610in
    \rlap{\kern 0.610in\lower\graphtemp\hbox to 0pt{\hss $\bullet$\hss}}%
    \graphtemp=.5ex
    \advance\graphtemp by 0.356in
    \rlap{\kern 0.356in\lower\graphtemp\hbox to 0pt{\hss $\bullet$\hss}}%
    \graphtemp=.5ex
    \advance\graphtemp by 0.102in
    \rlap{\kern 0.356in\lower\graphtemp\hbox to 0pt{\hss $\bullet$\hss}}%
    \special{pn 8}%
    \special{pa 356 1119}%
    \special{pa 356 864}%
    \special{pa 610 610}%
    \special{pa 356 356}%
    \special{pa 356 102}%
    \special{fp}%
    \graphtemp=.5ex
    \advance\graphtemp by 1.119in
    \rlap{\kern 1.119in\lower\graphtemp\hbox to 0pt{\hss $\bullet$\hss}}%
    \graphtemp=.5ex
    \advance\graphtemp by 0.864in
    \rlap{\kern 1.119in\lower\graphtemp\hbox to 0pt{\hss $\bullet$\hss}}%
    \graphtemp=.5ex
    \advance\graphtemp by 0.610in
    \rlap{\kern 0.864in\lower\graphtemp\hbox to 0pt{\hss $\bullet$\hss}}%
    \graphtemp=.5ex
    \advance\graphtemp by 0.356in
    \rlap{\kern 1.119in\lower\graphtemp\hbox to 0pt{\hss $\bullet$\hss}}%
    \graphtemp=.5ex
    \advance\graphtemp by 0.102in
    \rlap{\kern 1.119in\lower\graphtemp\hbox to 0pt{\hss $\bullet$\hss}}%
    \special{pa 1119 1119}%
    \special{pa 1119 864}%
    \special{pa 864 610}%
    \special{pa 1119 356}%
    \special{pa 1119 102}%
    \special{fp}%
    \graphtemp=.5ex
    \advance\graphtemp by 0.864in
    \rlap{\kern 0.102in\lower\graphtemp\hbox to 0pt{\hss $\bullet$\hss}}%
    \graphtemp=.5ex
    \advance\graphtemp by 0.356in
    \rlap{\kern 0.102in\lower\graphtemp\hbox to 0pt{\hss $\bullet$\hss}}%
    \graphtemp=.5ex
    \advance\graphtemp by 0.864in
    \rlap{\kern 1.373in\lower\graphtemp\hbox to 0pt{\hss $\bullet$\hss}}%
    \graphtemp=.5ex
    \advance\graphtemp by 0.356in
    \rlap{\kern 1.373in\lower\graphtemp\hbox to 0pt{\hss $\bullet$\hss}}%
    \special{pa 102 864}%
    \special{pa 356 864}%
    \special{fp}%
    \special{pa 102 356}%
    \special{pa 356 356}%
    \special{fp}%
    \special{pa 1119 864}%
    \special{pa 1373 864}%
    \special{fp}%
    \special{pa 1119 356}%
    \special{pa 1373 356}%
    \special{fp}%
    \special{pa 610 610}%
    \special{pa 864 610}%
    \special{fp}%
    \graphtemp=.5ex
    \advance\graphtemp by 0.508in
    \rlap{\kern 0.610in\lower\graphtemp\hbox to 0pt{\hss 0\hss}}%
    \graphtemp=.5ex
    \advance\graphtemp by 0.508in
    \rlap{\kern 0.864in\lower\graphtemp\hbox to 0pt{\hss 1\hss}}%
    \graphtemp=.5ex
    \advance\graphtemp by 0.254in
    \rlap{\kern 0.102in\lower\graphtemp\hbox to 0pt{\hss 0\hss}}%
    \graphtemp=.5ex
    \advance\graphtemp by 0.102in
    \rlap{\kern 0.458in\lower\graphtemp\hbox to 0pt{\hss 1\hss}}%
    \graphtemp=.5ex
    \advance\graphtemp by 0.763in
    \rlap{\kern 0.102in\lower\graphtemp\hbox to 0pt{\hss 1\hss}}%
    \graphtemp=.5ex
    \advance\graphtemp by 1.119in
    \rlap{\kern 0.458in\lower\graphtemp\hbox to 0pt{\hss 0\hss}}%
    \graphtemp=.5ex
    \advance\graphtemp by 0.763in
    \rlap{\kern 0.356in\lower\graphtemp\hbox to 0pt{\hss 0\hss}}%
    \graphtemp=.5ex
    \advance\graphtemp by 0.356in
    \rlap{\kern 0.458in\lower\graphtemp\hbox to 0pt{\hss 1\hss}}%
    \graphtemp=.5ex
    \advance\graphtemp by 0.763in
    \rlap{\kern 1.373in\lower\graphtemp\hbox to 0pt{\hss 0\hss}}%
    \graphtemp=.5ex
    \advance\graphtemp by 1.119in
    \rlap{\kern 1.220in\lower\graphtemp\hbox to 0pt{\hss 1\hss}}%
    \graphtemp=.5ex
    \advance\graphtemp by 0.254in
    \rlap{\kern 1.373in\lower\graphtemp\hbox to 0pt{\hss 0\hss}}%
    \graphtemp=.5ex
    \advance\graphtemp by 0.102in
    \rlap{\kern 1.220in\lower\graphtemp\hbox to 0pt{\hss 1\hss}}%
    \graphtemp=.5ex
    \advance\graphtemp by 0.284in
    \rlap{\kern 1.191in\lower\graphtemp\hbox to 0pt{\hss 0\hss}}%
    \graphtemp=.5ex
    \advance\graphtemp by 0.763in
    \rlap{\kern 1.119in\lower\graphtemp\hbox to 0pt{\hss 1\hss}}%
    \graphtemp=.5ex
    \advance\graphtemp by 0.864in
    \rlap{\kern 1.881in\lower\graphtemp\hbox to 0pt{\hss $\bullet$\hss}}%
    \graphtemp=.5ex
    \advance\graphtemp by 0.864in
    \rlap{\kern 2.136in\lower\graphtemp\hbox to 0pt{\hss $\bullet$\hss}}%
    \graphtemp=.5ex
    \advance\graphtemp by 1.119in
    \rlap{\kern 2.136in\lower\graphtemp\hbox to 0pt{\hss $\bullet$\hss}}%
    \special{pa 1881 864}%
    \special{pa 2136 864}%
    \special{pa 2136 1119}%
    \special{pa 1881 864}%
    \special{fp}%
    \graphtemp=.5ex
    \advance\graphtemp by 0.356in
    \rlap{\kern 1.881in\lower\graphtemp\hbox to 0pt{\hss $\bullet$\hss}}%
    \graphtemp=.5ex
    \advance\graphtemp by 0.102in
    \rlap{\kern 2.136in\lower\graphtemp\hbox to 0pt{\hss $\bullet$\hss}}%
    \graphtemp=.5ex
    \advance\graphtemp by 0.356in
    \rlap{\kern 2.136in\lower\graphtemp\hbox to 0pt{\hss $\bullet$\hss}}%
    \special{pa 1881 356}%
    \special{pa 2136 102}%
    \special{pa 2136 356}%
    \special{pa 1881 356}%
    \special{fp}%
    \graphtemp=.5ex
    \advance\graphtemp by 0.356in
    \rlap{\kern 2.644in\lower\graphtemp\hbox to 0pt{\hss $\bullet$\hss}}%
    \graphtemp=.5ex
    \advance\graphtemp by 0.102in
    \rlap{\kern 2.644in\lower\graphtemp\hbox to 0pt{\hss $\bullet$\hss}}%
    \graphtemp=.5ex
    \advance\graphtemp by 0.356in
    \rlap{\kern 2.898in\lower\graphtemp\hbox to 0pt{\hss $\bullet$\hss}}%
    \special{pa 2644 356}%
    \special{pa 2644 102}%
    \special{pa 2898 356}%
    \special{pa 2644 356}%
    \special{fp}%
    \graphtemp=.5ex
    \advance\graphtemp by 1.119in
    \rlap{\kern 2.644in\lower\graphtemp\hbox to 0pt{\hss $\bullet$\hss}}%
    \graphtemp=.5ex
    \advance\graphtemp by 0.864in
    \rlap{\kern 2.644in\lower\graphtemp\hbox to 0pt{\hss $\bullet$\hss}}%
    \graphtemp=.5ex
    \advance\graphtemp by 0.864in
    \rlap{\kern 2.898in\lower\graphtemp\hbox to 0pt{\hss $\bullet$\hss}}%
    \special{pa 2644 1119}%
    \special{pa 2644 864}%
    \special{pa 2898 864}%
    \special{pa 2644 1119}%
    \special{fp}%
    \graphtemp=.5ex
    \advance\graphtemp by 0.610in
    \rlap{\kern 2.390in\lower\graphtemp\hbox to 0pt{\hss $\bullet$\hss}}%
    \special{pa 2136 864}%
    \special{pa 2136 356}%
    \special{pa 2390 610}%
    \special{pa 2136 864}%
    \special{fp}%
    \special{pa 2644 864}%
    \special{pa 2644 356}%
    \special{pa 2390 610}%
    \special{pa 2644 864}%
    \special{fp}%
    \graphtemp=.5ex
    \advance\graphtemp by 0.508in
    \rlap{\kern 2.390in\lower\graphtemp\hbox to 0pt{\hss 0\hss}}%
    \graphtemp=.5ex
    \advance\graphtemp by 0.763in
    \rlap{\kern 1.881in\lower\graphtemp\hbox to 0pt{\hss 1\hss}}%
    \graphtemp=.5ex
    \advance\graphtemp by 0.936in
    \rlap{\kern 2.208in\lower\graphtemp\hbox to 0pt{\hss 0\hss}}%
    \graphtemp=.5ex
    \advance\graphtemp by 1.119in
    \rlap{\kern 2.237in\lower\graphtemp\hbox to 0pt{\hss 0\hss}}%
    \graphtemp=.5ex
    \advance\graphtemp by 0.254in
    \rlap{\kern 1.881in\lower\graphtemp\hbox to 0pt{\hss 0\hss}}%
    \graphtemp=.5ex
    \advance\graphtemp by 0.102in
    \rlap{\kern 2.237in\lower\graphtemp\hbox to 0pt{\hss 1\hss}}%
    \graphtemp=.5ex
    \advance\graphtemp by 0.284in
    \rlap{\kern 2.208in\lower\graphtemp\hbox to 0pt{\hss 1\hss}}%
    \graphtemp=.5ex
    \advance\graphtemp by 0.284in
    \rlap{\kern 2.572in\lower\graphtemp\hbox to 0pt{\hss 2\hss}}%
    \graphtemp=.5ex
    \advance\graphtemp by 0.102in
    \rlap{\kern 2.746in\lower\graphtemp\hbox to 0pt{\hss 0\hss}}%
    \graphtemp=.5ex
    \advance\graphtemp by 0.254in
    \rlap{\kern 2.898in\lower\graphtemp\hbox to 0pt{\hss 1\hss}}%
    \graphtemp=.5ex
    \advance\graphtemp by 1.119in
    \rlap{\kern 2.542in\lower\graphtemp\hbox to 0pt{\hss 2\hss}}%
    \graphtemp=.5ex
    \advance\graphtemp by 0.792in
    \rlap{\kern 2.716in\lower\graphtemp\hbox to 0pt{\hss 0\hss}}%
    \graphtemp=.5ex
    \advance\graphtemp by 0.763in
    \rlap{\kern 2.898in\lower\graphtemp\hbox to 0pt{\hss 1\hss}}%
    \graphtemp=.5ex
    \advance\graphtemp by 1.119in
    \rlap{\kern 0.712in\lower\graphtemp\hbox to 0pt{\hss T\hss}}%
    \graphtemp=.5ex
    \advance\graphtemp by 1.119in
    \rlap{\kern 3.000in\lower\graphtemp\hbox to 0pt{\hss L(T)\hss}}%
    \hbox{\vrule depth1.220in width0pt height 0pt}%
    \kern 3.000in
  }%
}%
}

\begin{proof}
($\Rightarrow$)  Suppose that for some $j$ we have $|P_j \cap U| < m_j$.  We construct a coloring of $V(T)$ in which the center receives one color.  Then by Proposition \ref{monocent}, $D(T) = D(L(T))$.  We first observe that for every $\phi \in \aut(T)$, if $s \in P_i$ then $\phi(s) \in P_i$ as well.  Now suppose that $\phi \in \aut(T)$ transposes $u$ and $w$.  Then for each $s \in U$, $\phi(s) \in W$ and for each $s \in W$, $\phi(s) \in U$.  Therefore, if $|P_j \cap U| \neq |P_j \cap W|$ then there are no automorphisms of $T$ that transpose the two center vertices.  In this case, we can simply start with any optimal SB coloring for $T$ and recolor the center vertices to the same color.  

So we may assume that $|P_j \cap U| = |P_j \cap W| < m_j$.  We begin with an optimal SB coloring $f$ of $T$ and we recolor the center to the same color.  Then we recolor the subgraphs in $P_j \cap W$.  Since $|P_j \cap U| < m_j$ there is a rooted SB class that is not represented by $f$ over $P_j \cap U$.  We recolor $P_j \cap W$ (if necessary) so that this rooted SB class appears.  As we saw earlier, any automorphism which transposes the center must map $P_j \cap U$ into $P_j \cap W$.  However, no such automorphism can agree with $f$ because $f$ restricts to an SB coloring of some subtree in $P_j \cap W$ which is essentially distinct from those SB colorings that appear on subtrees in $P_j \cap U$.  Therefore, we have an SB coloring in which the center receives one color.

($\Leftarrow$)  Now suppose that the hypotheses hold and that we have a $k$-coloring of $L(T)$.  Suppose that the coloring distinguishes $L(T)$.  Let $P'_i = \{S'_j : S_j \in P_i\}$ and $U',W'$ be defined similarly.  Each $S'_i$ must be colored so as to disagree with $A'_i$ and there are exactly $D(S'_i;k) = D(S_i;k)$ essentially distinct ways to do this.  However in $P'_i \cap U'$ we find exactly this many subgraphs and therefore each of these coloring types must appear.  The same goes for $P'_i \cap W'$, so we can associate each subgraph $S'_j \in P'_i \cap U'$ with a subgraph $S'_l \in P'_i \cap W'$ so that they have equivalent colorings.  Note also that if a coloring class is represented more than once in $P'_i$ then there is a color-preserving automorphism transposing the subgraphs whose colorings are in the same class and leaving the rest of the graph fixed.  Hence, the association is a one-one correspondance.  For each $s \in P'_i \cap U'$, and $t \in P'_i \cap W$ which is associated with $s$, let $\beta_s:V(s) \to V(t)$ be a color preserving isomorphism where $\beta_s(e_{v_s}) = e_{v_t}$.

Let $\tau(\{u,w\}) = \{u,w\}$ and let $\mathbb{B} = \{\beta_s : s \in U'\} \cup \{\tau\}$ and $\beta = \bigcup{\mathbb{B}}$, where the union denotes that $\beta$ restricts to each function in $\mathbb{B}$ on its domain.  We claim that $\beta$ is a color preserving automorphism of $L(T)$.  We see that $\beta$ is a well-defined bijection $V(L(T)) \to V(L(T))$ since the above association is bijective and $V(S'_i) \cap V(S'_j) = \emptyset$ when $i \neq j$.  Clearly, $\beta$ is color preserving.  Suppose that $v,v' \in V(L(T))$ are adjacent.  If $v$ and $v'$ are both in the same subgraph $s \in S'$ then $\beta(v)$ and $\beta(v')$ are adjacent since $\beta_s$ was an isomorphism.  Suppose that $v,v' \neq \{u,w\}$ and are not in the same $s \in S'$.  Then $v$ and $v'$ are both adjacent to $\{u,w\}$ and sit in subgraphs either both in $U'$ or both in $W'$.  Then $\beta(v)$ and $\beta(v')$ will also be adjacent to $\{u,w\}$ since it is fixed and will both be in $U'$ or in $W'$ and are therefore adjacent.  Finally if $v = \{u,w\}$ then $v' = v_j$ for some $j$ and therefore $\beta(v') = v_i$ for some $i$ and therefore $\beta(v)$ and $\beta(v')$ are adjacent.  So we obtain a color preserving automorphism from any $k$-coloring of $L(T)$.  We conclude that $D(T) < D(L(T))$.
\end{proof}

\section{Further comments and questions}
We saw that the iterated line graph is very rigid with respect to automorphisms.  All of the dynamics are at the level of a fixed collection of subgraphs changing places.  We can investigate this idea further.  For example, say that $G$ is a graph such that $\{P_i\}_{i=1}^k$ is a partition of $V(G)$ where every automorphism $\phi \in \aut(G)$ can be represented by a map $\pi \in S_k$ where $\phi(P_i) = P_{\pi(i)}$.  Suppose that $\aut(G)$ is isomorphic in this way to $H \leq S_k$.  Let $Tr = \{i: P_i \text{ has a nontrivial } H\text{-orbit}\}$ and $r = \min_{i \in Tr}{n(G_i)}$.  Let $m$ be the least number such that $\binom{r + m - 1}{m - 1} \geq k$.  Suppose that if $P_i$ has a trivial $H$-orbit then every automorphism restricts to the trivial automorphism on $P_i$.  Then we require at most $m$ colors to distinguish $G$ since there are at least $\binom{r + m - 1}{m - 1}$ ways to color the vertices of $P_j$ when $P_j$ has a nontrivial orbit.  Hence we can color each of the $P_i$ which have a nontrivial orbit with a different distribution of the $m$ colors and because $\aut(G)$ is represented faithfully in $S_k$, we obtain a distinguishing in the same vain as Theorem \ref{mainthm}.  Are there other reasons why the distinguishing number of a graph might be 2?

We obtain a bound on the number of iterations of the line graph required to reach distinguishing number 2, although it almost certainly is not optimal.  Does an optimal bound exist?  The pursuit of such a bound would likely proceed through the territory of short term behavior, an area where there are likely to be interesting phenomena.

We have seen one way that the distinguishing number of a graph can increase under $L$.  In computing examples, it seems that more often $D(G)$ is decreasing monotonically.  Can the graphs $G$ with the property that $D(L(G)) > D(G)$ be characterized nicely?  We have not been able to construct any additional graphs with this property.

Recall the importance of $D(G;k)$ as a parameter.  We used it in section four, but there is almost certainly more to it.  Given a simple graph $G$ and $v \in V(G)$ let $D_k(G)$ be the collection of SB colorings $f:G \to \{1,2\dotsc,k\}$ with respect to the stabilizer of $v$.  Say $f,g \in D_k(G)$ are equivalent, $f \sim g$, if there exists $\phi \in \aut(G)$ such that $f(v) = g(\phi(v))$ for all $v \in V(G)$.  Let $D(G;k)$ be the number of equivalence classes in $D_k(G)/\sim$.  What do $D_k(G)$, $D_k(G)/\sim$, and $D(G;k)$ look like?  Does $D(G;k)$ have any interesting properties?

\bibliographystyle{amsplain}
\bibliography{dnilg}

\def\soft#1{\leavevmode\setbox0=\hbox{h}\dimen7=\ht0\advance \dimen7
  by-1ex\relax\if t#1\relax\rlap{\raise.6\dimen7
  \hbox{\kern.3ex\char'47}}#1\relax\else\if T#1\relax
  \rlap{\raise.5\dimen7\hbox{\kern1.3ex\char'47}}#1\relax \else\if
  d#1\relax\rlap{\raise.5\dimen7\hbox{\kern.9ex \char'47}}#1\relax\else\if
  D#1\relax\rlap{\raise.5\dimen7 \hbox{\kern1.4ex\char'47}}#1\relax\else\if
  l#1\relax \rlap{\raise.5\dimen7\hbox{\kern.4ex\char'47}}#1\relax \else\if
  L#1\relax\rlap{\raise.5\dimen7\hbox{\kern.7ex
  \char'47}}#1\relax\else\message{accent \string\soft \space #1 not
  defined!}#1\relax\fi\fi\fi\fi\fi\fi}
\providecommand{\bysame}{\leavevmode\hbox to3em{\hrulefill}\thinspace}
\providecommand{\MR}{\relax\ifhmode\unskip\space\fi MR }
\providecommand{\MRhref}[2]{%
  \href{http://www.ams.org/mathscinet-getitem?mr=#1}{#2}
}
\providecommand{\href}[2]{#2}
\begin{thebibliography}{1}

\bibitem{Alb}
Michael~O. Albertson and Karen~L. Collins, \emph{Symmetry breaking in graphs},
  Electron. J. Combin. \textbf{3} (1996), no.~1, Research Paper 18, approx.\ 17
  pp.\ (electronic).

\bibitem{HH1}
Stephen~G. Hartke and Aparna~W. Higgins, \emph{Maximum degree growth of the
  iterated line graph}, Electron. J. Combin. \textbf{6} (1999), Research paper
  28, 9 pp.\ (electronic).

\bibitem{HH2}
\bysame, \emph{Minimum degree growth of the iterated line graph}, Ars Combin.
  \textbf{69} (2003), 275--283.

\bibitem{KnNi}
Martin Knor and {\soft{L}}udov{\'{\i}}t Niepel, \emph{Connectivity of iterated
  line graphs}, Discrete Appl. Math. \textbf{125} (2003), no.~2-3, 255--266.

\bibitem{KNS}
{\soft{L}}.~Niepel, M.~Knor, and L'. {\v{S}}olt{\'e}s, \emph{Distances in
  iterated line graphs}, Ars Combin. \textbf{43} (1996), 193--202.

\bibitem{Sab}
Gert Sabidussi, \emph{Graph derivatives}, Math. Z. \textbf{76} (1961),
  385--401.

\end{thebibliography}
\end{document}